\documentclass[12pt]{article}
\usepackage{latexsym,amsfonts}
\newtheorem{Theorem}{Theorem}[section]

\newtheorem{Lemma}[Theorem]{Lemma}
\newtheorem{Proposition}[Theorem]{Proposition}

\newtheorem{Remark}[Theorem]{Remark}
\newtheorem{Open Problem}[Theorem]{Open Problem}

\newcommand{\mysection}[1]{\section{#1}\setcounter{equation}{0}}
\def\pf{\noindent{\bf Proof. }}

\newcommand{\R}{\mathbb{R}}
\newcommand{\C}{\mathbb{C}}
\def\be{{}\begin{equation}}
\def\ee{{}\end{equation}}

\def\la{{}\langle}
\def\ra{{}\rangle}

\begin{document}

\title{
Orbital stability of standing wave solution for\\ a quasilinear
Schr\"{o}dinger equation}

\author{Boling Guo$^1$\ and\ Jianqing Chen$^{1,2}$
\thanks{Corresponding author: Institute of Applied Physics and
Computational Mathematics, PO Box 8009(28 Branch),  Beijing
100088, P. R. China. Supported in part by Tian Yuan Foundation of
Chinese NSF
(A0324612).}\\
 \small 1 Institute of Applied Physics and
Computational
Mathematics, PO Box 8009,\\
\small   Beijing 100088, P. R. China\\
\small 2 School of Mathematics and Computer Science, Fujian Normal
University,\\
 \small Fuzhou, 350007, P. R. China }

\date{}
\maketitle

\begin{minipage}{13cm}
{\small {\bf Abstract:}
 Via minimization arguments and Concentration Compactness Principle,
 we prove the orbital stability of standing wave solutions
 for a class of quasilinear Schr\"{o}dinger equation arising from physics.
}
\end{minipage}

\mysection{Introduction} This paper is motivated by the recent
interests on the following type of quasilinear Schr\"{o}dinger
equation
\begin{equation}\label{eq11}
i\partial_tz = -\Delta z + V(x)z -k(\Delta(|z|^2))z -
\theta|z|^{p-2}z,
\end{equation}
where {\it i} is the imaginary unit, $p > 2$ and $k,\ \theta\in
\R_+$, $\Delta z =\sum_{j=1}^N\partial^2z/\partial x_j^2$ the
standard Laplacian operator. $z:=z(x,t):\ \R^N\times\R_+\to \C$ is
a complex-valued function.

Problems of this kind arise naturally from various domains of
mathematical physics and have been derived as models of several
physical phenomena in the theory of superfluid film and in
dissipative quantum mechanics (see e.g.  Kurihura \cite{ku},
Nakamura \cite{na} ). For more physical motivations and more
references dealing with applications, we refer the interested
readers to Lange et al \cite{lpt}, Poppenberg et al \cite{psw} and
the references therein.

Via critical point theory, finding a standing wave of the form
$z(x,t)=e^{i\mu t}u(x)$ of problem (\ref{eq11}) is equivalent to
solve the following elliptic equation
\begin{equation}\label{eq12}
-\Delta u + (V(x) + \mu)u -k(\Delta(|u|^2))u - \theta|u|^{p-2}u=0.
\end{equation}
Under the basic assumption of $\inf_{x\in\R^N}V(x)+\mu > 0$, the
existence of nontrivial solutions of (\ref{eq12}) has been studied
in the last two decades. Some of them are due to Strauss
\cite{str} and Rabinowitz \cite{rab} in the case of $k=0$ and
Ambrosetti-Liu-Wang-Wang \cite{aw, lw, lww1, lww2} in the case of
$k\neq 0$ with some additional assumptions on $V(x)$ and $p$.
Using the basic function space $H^1(\R^N)=W^{1,2}(\R^N)$ and
setting
$$X=\{u\in H^1(\R^N);\ \int_{\R^N}|u|^2|\nabla u|^2dx < \infty,\ \int V(x)|u|^2<\infty\},$$
we can define functionals:
$$F_1(u)={1\over 2}\int(|\nabla u|^2+V(x)|u|^2),\quad\quad F_2(u)={1\over
2}\int|u|^2,$$
$$F_3(u)={1\over 4}\int|\nabla|u|^2|^2,\quad\quad F_4(u)={1\over
p}\int|u|^p,$$
$$F_5(u)=F_1(u)+\mu F_2(u)+kF_3(u)-\theta F_4(u),$$
$$F_6(u)=\int[|\nabla u|^2+(V(x)+\mu)|u|^2+k|\nabla|u|^2|^2-\theta|u|^p].$$
Problem (\ref{eq12}) has a formal variational structure. For any
$\phi\in\mathcal{D}(\R^N)$, $F_j$, $j=1,\ \cdots,\ 6$, have
directional derivatives at $u$ in the direction $\phi$, denoted by
$\la F'_j(u),\phi\ra$. We say that $u$ is a weak solution of
(\ref{eq12}) if and only if for any $\phi\in\mathcal{D}(\R^N)$ and
some $\mu,\ \theta$, there holds
\begin{equation}\label{eq13}
\la F'_1(u) + \mu F'_2(u) + kF_3'(u)-\theta F_4'(u),\phi\ra = 0.
\end{equation}
Up to now, we have three methods to
study the existence of standing wave solutions for problems of
this kind.\\

\noindent{\it Method I:} Study the minimization problem
\begin{equation}\label{eq14}
m_I=\inf\{F_1(u)+\mu F_2(u)+kF_3(u);\ u\in X,\
F_4(u)=\lambda_I>0\}.
\end{equation}
If $m_I$ is achieved, then the minimizer is  a
solution of (\ref{eq12}) for some $\theta$.\\

\noindent{\it Method II:} Consider the minimization problem
\begin{equation}\label{eq15}
m_{II}=\inf\{F_1(u)+kF_3(u)-\theta F_4(u);\ u\in X,\
F_2(u)=\lambda_{II}>0\}.
\end{equation}
If $m_{II}$ is achieved, then again the minimizer is  a
solution of (\ref{eq12}) for some $\mu$.\\

\noindent{\it Method III:} Study the minimization problem
\begin{equation}\label{eq16}
m_{III}=\inf\{F_5(u);\ 0\neq u\in X,\ F_6(u)=0\}.
\end{equation}
If $m_{III}$ is achieved, then the minimizer is a solution of
(\ref{eq12}).\\

When $k=0$, the three methods mentioned above are almost equivalent
in the sense that they can be changed from one to another by
scaling. However, when $k\neq 0$, the term introduced by $F_3(u)$
make the problem substantially different from the case of $k=0$
since $F_3$ is non-convex and the scaling argument does not work.
{\it Method III} has been used in \cite{lw, lww1, lww2} to study the
existence of solutions of (\ref{eq12}). But it seems that it is too
difficult to be used to study orbital stability (see precise
definition in Section 3) of standing waves of (\ref{eq11}) since
$F_6$ usually does not satisfy any conservative laws (see Lemma
\ref{le21}). Hence our purpose here is to take {\it Method II} to
study the existence and orbital stability of standing waves of
(\ref{eq11}) under the basic assumption $\inf_{x\in\R^N}V(x) \geq
0$.

We recall that in the case of $k=0$, the existence and orbital
stability of standing wave solution of (\ref{eq11}) has been
studied by Cazenave and Lions \cite{cal} via concentration
compactness principle. However, the presence of the quasilinear
term $(\Delta(|z|^2))z$ makes the problem substantially different
from the semilinear case. For example, when $k=0$, one can use
$\{u;\ \int|u|^2=\hbox{const.}\}$ or $\{u;\
\int|u|^p=\hbox{const.}\}$ as constraint to study a corresponding
minimization problem. Then the study of orbital stability can be
proceeded after standard scaling argument. But the term introduced
by $(\Delta(|z|^2))z$ makes the scaling arguments fail. Hence we
need to restrict the arguments to a function space which consists
of functions with radial symmetry. The main results are Theorem
\ref{th32} and Theorem \ref{th33}, which are contained in Section
3 after some preliminaries given in Section 2. In Section 4, we
study problem (\ref{eq11}) in one spatial dimension.
 Some improvements of Theorem \ref{th32} and Theorem \ref{th33} are given, see
 Theorem \ref{th46} and Theorem \ref{th48}.\\

\mysection{Preliminaries} Throughout this paper, all integrals are
taken over $\R^N$ (or $\R$ which were understood from the
contexts) unless stated otherwise. All $dx$ in the integrals are
omitted. $o(1)$ will denote a generic infinitesimal as
$n\to\infty$. $\rightarrow$ denotes the strong convergence and
$\rightharpoonup$ the weak convergence. $B(x,R)$ will denote a
ball centered at $x$ with radials $R$ and $B(0,R)$ is simply
denoted by $B_R$. To continue, we list some lemmas which will be
useful in what follows. First
\begin{Lemma}[Conservative laws]\label{le21} Let $z(x,t)$ be a
solution of (\ref{eq11}) with initial value $z_0$. Then there hold
\begin{equation}\label{eq21}
\int|z(x,t)|^2 = \int|z_0|^2\equiv \hbox{const.}\quad
\hbox{Conservation of ``mass"}
\end{equation}
\begin{equation}\label{eq22}
\begin{array}{rl}&\displaystyle\int[{1\over 2}(|\nabla z|^2+V|z|^2)+
{k\over 4}|\nabla|z|^2|^2-{\theta\over p}|z|^p]\\
 &= \displaystyle\int[{1\over 2}(|\nabla
z_0|^2+V|z_0|^2)+ {k\over 4}|\nabla|z_0|^2|^2-{\theta\over p}|z_0|^p]\\
& \equiv \hbox{const.}\quad \hbox{Conservation of
``energy"}\end{array}
\end{equation}
\end{Lemma}
\pf This kind of result should be known, see e.g. \cite{lpt}.$\Box$
\begin{Lemma}\label{le22}\cite{lio}
Let $(\rho_n)$ be bounded in $L^1(\R)$. By extracting a
subsequence, we may assume that $(\rho_n)$ satisfies one of the
following two possibilities:
\begin{description}
    \item[(i)] (Vanishing)
    $\lim_{n\to\infty}\sup_{y\in\R}\displaystyle\int_{y+B_R}\rho_n(x)dx=0$
    for all $0<R<+\infty$;
    \item[(ii)] (Non-vanishing) There exist $\alpha > 0$, $R <
    +\infty$ and $(y_n)\subset\R$ such that
    $$\lim_{n\to\infty}\displaystyle\int_{y_n+B_R}\rho_n(x)dx\geq \alpha
    >0.$$
\end{description}
\end{Lemma}
\begin{Lemma}\label{le23}\cite{lio}
Let $(u_n)$ be bounded in $H^1(\R)$. Assume that for some $q > 2$
and $R > 0$,
$$\sup_{y\in\R}\displaystyle\int_{y+B_R}|u_n|^qdx \to 0\quad \hbox{as}\quad n\to \infty,$$
then $u_n\to 0$ in $L^\beta(\R)$ for any $\beta \geq 2$.
\end{Lemma}
\begin{Lemma}\cite[Lemma 2 in Page 333]{psw}\label{le24}
If $u_n\rightharpoonup u$ in $H^1(\R)$ and $u_n\to u$ a.e. in
$\R$, then
$$\liminf_{n\to\infty}\int|(|u_n|^2)'|^2 \geq
 \int|(|u|^2)'|^2 + \liminf_{n\to\infty}\int|(|u_n-u|^2)'|^2.$$
\end{Lemma}

Denote $H^1_r(\R^N)=\{u\in H^1(\R^N);\ u(x)=u(|x|)\}$ and $H=\{u\in
H^1_r(\R^N);\ \int V(x)|u|^2<+\infty\}$ with the norm
$\|u\|^2=\int[|\nabla u|^2+V(x)|u|^2+|u|^2]$. Then the following
lemma is by now well known.
\begin{Lemma}\label{le25}
Let $N\geq 2$. Then the following embedding is compact,
$$H\hookrightarrow L^q(\R^N),\quad 2\leq q < 2^\ast,$$
where $2^\ast={{2N}\over{N-2}}$ for $N\geq 3$ and $\infty$ for
$N=2$.
\end{Lemma}
\begin{Lemma}\cite[Lemma 13 in Page 340]{psw}\label{le26}
If $(u_n)\subset H$ are such that $u_n\rightharpoonup u$ in
$H^1(\R^N)$, then
$$\liminf_{n\to\infty}\int|\nabla|u_n|^2|^2 \geq
 \int|\nabla|u|^2|^2.$$
\end{Lemma}

\mysection{The case of $N\geq 2$} In this section, we assume
$V(x)=V(|x|)$. We will
 first follow the line of {\it method II} to study the existence of standing
 wave of (\ref{eq11}) in the case of $\theta > 0$.
Then we will study the orbital stability of the standing wave. Now
for any $\lambda > 0$ fixed, consider the following minimization
problem
\begin{equation}\label{eq31}
m_r = \inf\{E(u)=F_1(u)+kF_3(u)-\theta F_4(u);\ u\in H,\
F_2(u)=\lambda\}.
\end{equation}
\begin{Lemma}\label{le31}
Let $(u_n)\subset H$ be a minimizing sequence of $m_r$. If $2 < p
< 2+{4\over N}$, then $(u_n)$ is bounded in $H$.
\end{Lemma}
\pf Since $(u_n)\subset H$ is a minimizing sequence of $m_r$,
i.e.,
\begin{equation}\label{eq32}
m_r+o(1)=F_1(u_n)+kF_3(u_n)-\theta F_4(u_n);\quad
F_2(u_n)=\lambda.
\end{equation}
By Sobolev and interpolation inequalities, we deduce that for any
$u\in H$ and $s=({1\over p}-{1\over{2^\ast}})/({1\over
2}-{1\over{2^\ast}})$,
$$
|u|_p^p \leq |u|_2^{ps}|u|_{2^\ast}^{p(1-s)}\leq
M|u|_2^{ps}\|u\|^{p(1-s)}.
$$ Young inequality implies that
$$
|u|_p^p\leq \varepsilon\|u\|^2 +
A_\varepsilon|u|_2^{2ps/(2-p+ps)}.
$$
In here we use the assumption $2 < p < 2+{4\over N}$. It follows
that
\begin{equation}\label{eq33}
m_r+1\geq F_1(u_n)+{k\over 4}\int|\nabla|u_n|^2|^2 -
{\varepsilon\over p}\|u_n\|^2 - {{A_\varepsilon}\over
p}|u_n|_2^{2ps/(2-p+ps)}.
\end{equation}
We obtain from choosing $\varepsilon$ small enough and using
$|u_n|_2^2=2\lambda$ that $m_r > -\infty$ and $(u_n)$ is bounded
in $H$.$\Box$

\begin{Theorem}\label{th32}
Suppose $N\geq 2$, $2 < p < 2 + {4\over N}$ and $k>0$. Then for any
fixed $\theta>0$ and $\lambda > 0$,  $m_r$ is achieved at some
$u_0$, i.e., $u_0(x)=u_0(|x|)$, $F_2(u_0)=\lambda$,
$F_1(u_0)+kF_3(u_0)-\theta F_4(u_0)=m_r$.
\end{Theorem}
\pf Let $(u_n)\subset H$ be a minimizing sequence of $m_r$. Lemma
\ref{le31} implies that $(u_n)$ is bounded in $H$. Going if
necessary to a subsequence, we can assume that $u_n\rightharpoonup
u_0$ in $H$ and $u_n\rightarrow u_0$ a.e. in $\R^N$. Brezis-Lieb
Lemma \cite{bl} implies that
\begin{equation}\label{eq34}
F_1(u_n) = F_1(u_n-u_0) + F_1(u_0) + o(1).
\end{equation}
Combining this with (\ref{eq31}), Lemma \ref{le25} and Lemma
\ref{le26} we have that
$$\begin{array}{rl}
m_r+o(1) &= F_1(u_n) + {k\over 4}\int|\nabla|u_n|^2|^2
-{\theta\over p}|u_n|_p^p\\
&\geq F_1(u_0)+kF_3(u_0)-\theta
F_4(u_0) + \lim_{n\to\infty}F_1(u_n-u_0)\\
&\geq F_1(u_0)+kF_3(u_0)-\theta F_4(u_0) \geq m_r.\end{array}$$ It
follows that $m_r=F_1(u_0)+kF_3(u_0)-\theta F_4(u_0)$ and
$\lim_{n\to\infty}F_1(u_n-u_0) = 0$. Thus we obtain from
$|u_n-u_0|_2^2\to 0$ that $\|u_n-u_0\|\to 0$, i.e., $u_n\to u_0$
in $H$. $\Box$

 The method of studying minimization problem (\ref{eq31}) is
not new, but it has the merit of studying orbital stability of
standing wave of (\ref{eq11}) in the spirit of Cazenave and Lions
\cite{cal}, see also Albert \cite{alb}. Roughly speaking, a set of
solutions of (\ref{eq11}) is said to be stable if any solution of
(\ref{eq11}) remains near the set whenever it starts near the set.
Now for any $\lambda
> 0$, we know from Theorem \ref{th32}
that the set of minimizers, denoted by $S_\lambda$, of the
minimization problem (\ref{eq31}) is not empty. Then for any $u_0\in
S_\lambda$, $\int|\nabla|u_0|^2|^2$ is finite. Hence the integral
$\int\nabla |u_0|^2\nabla(u_0v)$ exists for any
$v\in\mathcal{D}(\R^N)$ because
\begin{equation}\label{eq35}
\int|u_0||\nabla u_0|^2 + |u_0|^2|\nabla u_0|\leq
\int(1+|u_0|^2)|\nabla u_0|^2 + |u_0|^2(1+|\nabla u_0|^2) <
+\infty.
\end{equation}
Therefore, the standard proof of the Ljusternik's Theorem on
Lagrange multipliers \cite{bro} implies that there exist $\gamma$
such that $u_0$ is a weak solution of (\ref{eq12}) for
$\mu=-\gamma$. It follows that $z(x,t)=e^{i\mu t}u_0(x)$ is a
standing wave of (\ref{eq11}). Thus $e^{i\mu t}u_0(\cdot)$ is the
orbit of $u_0$. Moreover, for any $t\geq 0$, if $u\in S_\lambda$,
then $e^{i\mu t}u(x)\in S_\lambda$. Our orbital stability result
is
\begin{Theorem}\label{th33}
Assume that $N\geq 2$ and $2 < p < 2+{4\over N}$. Then for any
$\varepsilon > 0$, there exists $\delta
> 0$ such that if
$$\inf_{g\in S_\lambda}\|z_0 - g\| < \delta,$$ then the solution
$z(x,t)$ of (\ref{eq11}) with initial value $z(x,0)=z_0$ satisfies
$$\inf_{g\in S_\lambda}\|z(\cdot,t) - g\| < \varepsilon$$
for any $t\in[0, T^*)$ with $T^*<\infty$ or $T^*=\infty$.
\end{Theorem}
\pf Suppose the conclusion to be false. Then there exist a number
$\varepsilon_0 > 0$, a sequence $(\psi_n)$ of functions in $H$ and
a sequence of times $(t_n)$ such that
$$\inf_{g\in S_\lambda}\|\psi_n - g\| < {1\over n}$$
and $$\inf_{g\in S_\lambda}\|z_n(\cdot,t_n) - g\| \geq
\varepsilon_0$$ for all $n$, where $z_n(x,t_n)$ solves
(\ref{eq11}) with $z_n(x,0) = \psi_n$. Then since $\psi_n\to
 S_\lambda$ in $H$ and $m_r = E(g)$ for all
$g\in S_\lambda$ and $F_2(g)=\lambda$, we have $E(\psi_n)\to m_r$
and $F_2(\psi_n)\to \lambda$ as $n\to \infty$. Thus we can find a
sequence $\beta_n\to 1$ such that $F_2(\beta_n\psi_n) = \lambda$
for all $n$. It follows from Lemma \ref{le21} that the sequence
$q_n = \beta_nz_n(\cdot,t_n)$ satisfies $F_2(q_n) = \lambda$ and
$$\lim_{n\to \infty}E(q_n) = \lim_{n\to \infty}E(z_n(\cdot,t_n))
= \lim_{n\to \infty}E(\psi_n) = m_r$$ and is therefore a
minimizing sequence for $m_r$. The proof of Theorem \ref{th32}
implies that $q_n\to q_0$ strongly in $H$. Hence
$|q_0|_2^2=2\lambda$ and we have a sequence $(g_n)\subset
S_\lambda$ such that $\|q_n - g_n\| < {\varepsilon_0\over 2}$ for
$n$ large. But
$$\begin{array}{rl}
\varepsilon_0& \leq \|z_n(\cdot,t_n) - g_n\|\leq \|z_n(\cdot,t_n)
- q_n\|+\|q_n - g_n\|\\
&\leq |1-\beta_n|\|z_n(\cdot,t_n)\| +\|q_n - g_n\|<
{\varepsilon_0\over 2},\end{array}$$ which is a contradiction.$\Box$
\begin{Remark}\label{re34}
Please note that in the stating of this Theorem, we have tacitly
assumed that the potential $V(x),\ p$ and the initial data $z_0$
belong to a class in which unique solution of the initial-value
problem for (\ref{eq11}) exist for all $t\in[0, T^*)$ and some
$T^*<\infty$ or $T^*=\infty$. In view of local well-posedness
results for this kind of problems studied by Lange-Poppenberg
\cite{lpt} and the arguments developed by Cazenave \cite{caz}, we
will continue to make this assumptions without any comments.
Moreover, if $z_0$ is radial symmetric with respect to $x$, then
so is $z(x,t)$.
\end{Remark}
\begin{Remark}\label{le35}
Note that in Theorem \ref{th32} and Theorem \ref{th33}, $V(x)$ can
cover the case of $V(x)=|x|^2$. Problem (\ref{eq11}) in the case
of $V(x)=|x|^2$ and $k=0$ describes the Bose-Einstain condensate
with attractive interparticle interactions under magnetic trap as
the Gross-Pitaevski equation with a harmonic potential term
\cite{tsw}. Hence Theorem \ref{th33} generalizes \cite[Theorem
3.2]{zha}.
\end{Remark}
\begin{Remark}\label{le36}
As pointed out in \cite[page 4]{alb}, although the Concentration
Compactness method for proving the orbital stability of standing
wave has the advantage of requiring less analysis than others (see
e.g., \cite{gss} and the references therein), it produce a weak
result in that it only demonstrates stability of a set of
minimizing solutions without providing information on the
structure of the set, or distinguishing among its different
members. This kind of possible confusion will be clarified in the
case of $N=1$.
\end{Remark}

\mysection{The case of $N=1$}
 In this section, we will give a stronger orbital
stability of standing wave of problem (\ref{eq11}) in the case of
$V(x)\equiv 0$ and $N=1$, $k\geq 0$, $\theta>0$. At this time, we
write (\ref{eq11}) as
\begin{equation}\label{eq41}
i\partial_tz = -z'' - k(|z|^2)''z - \theta|z|^{p-2}z
\end{equation}
and denote by $''$ (resp. $'$) the second (resp. first) order
spatial derivatives. Finding a standing wave of (\ref{eq41}) of
the form $e^{i\mu t}u(x)$ is equivalent to solve the associated
elliptic problem
\begin{equation}\label{eq42}
-u'' + \mu u -k((|u|^2)'')u - \theta|u|^{p-2}u=0.
\end{equation}
Using the basic function space $H^1(\R)=W^{1,2}(\R)$ with the
standard norm $\|u\|^2=\int[|\nabla u|^2+|u|^2]$ and the
continuous embedding $H^1(\R)\hookrightarrow L^\infty(\R)$, we
know that the functional
\begin{equation}\label{eq43}
I(u)=\int[{1\over 2}|u'|^2+{k\over 4}|(|u|^2)'|^2-{\theta\over
p}|u|^p]
\end{equation}
is well defined on $H^1(\R)$. What's more, we have
\begin{Lemma}\cite[Lemma 1]{psw}\label{le41}
$I$ is of class $C^1$ on $H^1(\R)$.$\Box$
\end{Lemma}
For any $\lambda > 0$ fixed, we define
\begin{equation}\label{eq44}
\mathcal{M}=\{u\in H^1(\R);\ {1\over
2}\int_{\R}|u|^2=\lambda\},\quad\quad
m=\inf_{u\in\mathcal{M}}I(u).
\end{equation}
\begin{Lemma}\label{le42}
If $2 < p < 6$, then $-\infty < m < 0$.
\end{Lemma}
\pf Suppose $\psi(x)\in\mathcal{M}$. Then so is $\xi^{1\over
2}\psi(\xi x)$ for $\xi > 0$. Hence
$$\begin{array}{rl}
I(\xi^{1\over 2}\psi(\xi x))&=\displaystyle\int_{\R}[{1\over
2}\xi^3|\psi'(\xi x)|^2+{k\over 4}\xi^4|(|\psi(\xi
x)|^2)'|^2-{\theta\over
p}\xi^{p\over 2}|\psi(\xi x)|^p]\\
&=\displaystyle\int[{1\over 2}\xi^2|\psi'|^2+{k\over
4}\xi^3|(|\psi|^2)'|^2-{\theta\over p}\xi^{{p\over
2}-1}|\psi|^p].\end{array}$$ It follows from $2< p < 6$ and $m\leq
I(\xi^{1\over 2}\psi(\xi x))$ that $m < 0$ provided $\xi > 0$ small
enough. Now for any $u\in\mathcal{M}$, we can use Sobolev and
interpolation inequalities that for some $0 < \alpha < 1$,
\begin{equation}\label{eq45}
|u|_p^p\leq |u|_2^{p\alpha}|u|_\infty^{p(1-\alpha)}\leq
d|u|_2^{p\alpha}\|u\|^{p(1-\alpha)}.
\end{equation}
Young inequality implies that
\begin{equation}\label{eq46}
|u|_p^p=\varepsilon\|u\|^2+A_\varepsilon|u|_2^{2p\alpha/(2-p+p\alpha)}.
\end{equation}
In here, we use the assumption $2 < p < 6$. By choosing some
$u\in\mathcal{M}$ such that
\begin{equation}\label{eq47}
m+1\geq {1\over 2}\|u\|^2-\lambda+{k\over
4}\int_{\R}|(|u|^2)'|^2-{\theta\over
p}\varepsilon\|u\|^2-{\theta\over
p}A_\varepsilon|u|_2^{2p\alpha/(2-p+p\alpha)}
\end{equation}
and $\varepsilon$ small enough, we can get that $m >
-\infty$.$\Box$
\begin{Lemma}\label{le43}
Let $\{u_n\}$ be a minimizing sequence of $m$. Then there exist
$w\neq 0$ and $(y_n)\subset\R$ such that
$u_n(\cdot+y_n)\rightharpoonup w$ in $H^1(\R)$.
\end{Lemma}
\pf Let $\{u_n\}$ be a minimizing sequence of $m$. Then from the
proof of the second part of Lemma \ref{le42}, we know that $(u_n)$
is bounded in $H^1(\R)$. Using the fact that ${1\over
2}\int_\R|u_n|^2=\lambda > 0$ and Lemma \ref{le23}, we know that the
vanishing case does not occur for $\rho_n(x)=|u_n(x)|^2$. Therefore,
there exist $\varepsilon > 0$ and $R>0$ such that
\begin{equation}\label{eq48}
\liminf_{n\to\infty}\sup_{y\in\R}\int_{y-R}^{y+R}|u_n(x)|^2 \geq
\varepsilon > 0.
\end{equation}
We may assume that there are $(y_n)\subset\R$ and $R>0$ such that
\begin{equation}\label{eq49}
\liminf_{n\to\infty}\int_{y_n-R}^{y_n+R}|u_n(x)|^2 \geq
{\varepsilon\over 2} > 0.
\end{equation}
Hence $w_n(x)=u_n(x+y_n)$ satisfies $I(w_n)=I(u_n)$ and ${1\over
2}\int_{\R}|w_n|^2=\lambda$, i.e., the sequence $(w_n)$ is also a
minimizing sequence which satisfies
\begin{equation}\label{eq410}
\liminf_{n\to\infty}\int_{-R}^R|w_n(x)|^2 \geq {\varepsilon\over
2} > 0.
\end{equation}
$(w_n)$ is also bounded in $H^1(\R)$. Going if necessary to a
subsequence, we have $w_n\rightharpoonup w$ in $H^1(\R)$ and thus
$w_n\to w$ in $L^2_{loc}(\R)$. (\ref{eq410}) implies that $w\neq
0$.$\Box$
\begin{Lemma}\label{le44}
Let $w_n\rightharpoonup w$ in $H^1(\R)$ and $w_n\to w$ a.e. in
$\R$. Then
\begin{equation}\label{eq411}
\lim_{n\to\infty}(|w_n|_q^q - |w_n-w|_q^q) = |w|_q^q,\quad q\geq
2;
\end{equation}
\begin{equation}\label{eq412}
\lim_{n\to\infty}(|w_n'|_2^2 - |(w_n-w)'|_2^2) = |w'|_2^2.
\end{equation}
\end{Lemma}
\pf These are direct consequence of Brezis-Lieb Lemma
\cite{bl}.$\Box$
\begin{Proposition}\label{pr45}
Suppose $4\leq p < 6$. Then $m$ is achieved by some nonnegative
$w$, i.e., there is $0\neq w\in \mathcal{M}$ with $I(w)=m$ and
$w\geq 0$.
\end{Proposition}
\pf Let $(u_n)$ be a minimizing sequence of $m$. Lemma \ref{le43}
implies that there are $(y_n)$ and $w\neq 0$ such that
$w_n(x)=u_n(x+y_n)$ and $w_n\rightharpoonup w\neq 0$ in $H^1(\R)$.
To conclude the proof, we have to show that $w\in\mathcal{M}$,
i.e., ${1\over 2}|w|_2^2=\lambda$. Arguing by a contradiction, we
may assume $0 < {1\over 2}|w|_2^2 < \lambda$. For $v_n=w_n-w$, we
may assume from Lemma \ref{le44} that ${1\over 2}|v_n|_2^2 <
\lambda$. Denote $\tilde{w}=w/({1\over{2\lambda}})^{1\over
2}|w|_2$ and $\tilde{v}_n=v_n/({1\over{2\lambda}})^{1\over
2}|v_n|_2$. We have that
\begin{equation}\label{eq413}
\begin{array}{rl}
m&\geq \displaystyle\int_{\R}[{1\over
2}|w_n'|^2+{k\over 4}|(|w_n|^2)'|^2-{\theta\over p}|w_n|^p]\\
&\geq \displaystyle\int_{\R}[{1\over
2}|w'|^2+{k\over 4}|(|w|^2)'|^2-{\theta\over p}|w|^p]\\
&+\displaystyle\int_{\R}[{1\over 2}|v_n'|^2-{\theta\over
p}|v_n|^p]
+{k\over 4}\liminf_{n\to\infty}\displaystyle\int_{\R}|(|v_n|^2)'|^2\\
&=(\sqrt{1\over {2\lambda}}|w|_2)^2{\displaystyle\int_{\R}{1\over
2}|\tilde{w}'|^2}+{k\over 4}(\sqrt{1\over{2\lambda}}|w|_2)^4
{\displaystyle\int_{\R}|(|\tilde{w}|^2)'|^2} -(\sqrt{1\over
{2\lambda}}|w|_2)^p{\theta\over
p}{\displaystyle\int_{\R}|\tilde{w}|^p}\\
&+(\sqrt{1\over {2\lambda}}|v_n|_2)^2{\displaystyle\int_{\R}{1\over
2}|\tilde{v}_n'|^2}-(\sqrt{1\over {2\lambda}}|v_n|_2)^p{\theta\over
p}{\displaystyle\int_{\R}|\tilde{v}_n|^p}+{k\over
4}{\displaystyle\liminf_{n\to\infty}}(\sqrt{1\over
{2\lambda}}|v_n|_2)^4{\displaystyle\int_{\R}|(|\tilde{v}_n|^2)'|^2}\\
&>(\sqrt{1\over {2\lambda}}|w|_2)^4{\displaystyle\int_{\R}[{1\over
2}|\tilde{w}'|^2+{k\over 4}|(|\tilde{w}|^2)'|^2-{\theta\over
p}|\tilde{w}|^p]}\quad (\ \hbox{since}\ p\geq 4)\\
&+{\displaystyle\liminf_{n\to\infty}}(\sqrt{1\over
{2\lambda}}|v_n|_2)^4{\displaystyle\int_{\R}[{1\over
2}|\tilde{v}_n'|^2+{k\over 4}|(|\tilde{v}_n|^2)'|^2-{\theta\over
p}|\tilde{v}_n|^p]}\\
&=(\sqrt{1\over
{2\lambda}}|w|_2)^4I(w)+{\displaystyle\liminf_{n\to\infty}}(\sqrt{1\over
{2\lambda}}|v_n|_2)^4I(\tilde{v}_n)\\
&\geq m[(\sqrt{1\over
{2\lambda}}|w|_2)^4+{\displaystyle\liminf_{n\to\infty}}(\sqrt{1\over
{2\lambda}}|v_n|_2)^4]\\
&\geq m[(\sqrt{1\over
{2\lambda}}|w|_2)^2+{\displaystyle\liminf_{n\to\infty}}(\sqrt{1\over
{2\lambda}}|v_n|_2)^2]\quad (\ \hbox{since}\ m < 0)\\
&= m,
\end{array}
\end{equation}
which is a contradiction. In here the last equality follows from
using Lemma \ref{le44} and ${1\over 2}|w_n|_2^2=\lambda$. Up to
now we have proved that ${1\over 2}|w|_2^2=\lambda$, i.e., $w_n\to
w$ in $L^2(\R)$. Since $(w_n)$ is bounded in $H^1(\R)$, the
standard interpolation and Sobolev inequalities imply that $w_n\to
w$ in $L^p(\R)$. We now obtain from $(w_n)$ being a minimizing
sequence and Lemma \ref{le24} that
\begin{equation}\label{eq414}
\begin{array}{rl}
m&\geq \displaystyle\int_{\R}[{1\over
2}|w_n'|^2+{k\over 4}|(|w_n|^2)'|^2-{\theta\over p}|w_n|^p]\\
&\geq \displaystyle\int_{\R}[{1\over 2}|w'|^2+{k\over
4}|(|w|^2)'|^2-{\theta\over
p}|w|^p]+\displaystyle\lim_{n\to\infty}\int_{\R}{1\over
2}|(w_n-w)'|^2\\
&\geq m + \displaystyle\lim_{n\to\infty}\int_{\R}{1\over
2}|(w_n-w)'|^2.\end{array}\end{equation} Hence $\int_{\R}
|(w_n-w)'|^2\to 0$. Combining it with $w_n\to w$ in $L^2(\R)$ we
get that $w_n\to w$ in $H^1(\R)$. Note that if $(w_n)$ is a
minimizing sequence of $m$, then so is $(|w_n|)$. Therefore the
conclusion follows.$\Box$
\begin{Theorem}\label{th46}
Assume $4\leq p < 6$. Then there exists $\gamma < 0$ such that
(\ref{eq42}) admits a positive solution $u\in H_0^1(\R)$ for any
$\theta > 0$ and $\mu=-\gamma$.
\end{Theorem}
\pf With the help of Proposition \ref{pr45} and an argument in
front of Theorem \ref{th33}, we know that there exists a
nonnegative $u\in H^1(\R)$ with $|u|_2^2=2\lambda$ and a Lagrange
multiplier $\gamma\in\R$ such that
\begin{equation}\label{eq415}
-\int u''v- k\int(u^2)''uv - \theta|u|^{p-2}uv=\gamma\int uv,
\end{equation}
for all $v\in H^1(\R)$. Putting $v=u$ in (\ref{eq415}), we get
that
$$\gamma|u|_2^2=|u'|^2_2+k\int|(|u|^2)'|^2-\theta|u|_p^p.$$
It follows from $0 > m={1\over 2}|u'|_2^2+{k\over
4}\int|(|u|^2)'|^2-{\theta\over p}|u|_p^p$ and $p\geq 4$ that
$\gamma<0$. Hence $u\in H_0^1(\R)$ is a nonnegative solution of
(\ref{eq42}) for any $\mu=-\gamma$. The proof of positivity of $u$
is the same as \cite[Theorem 1]{psw} and is
omitted here. The proof is completed.$\Box$\\

As those pointed out in Remark \ref{re34}, we need to use
the following basic assumption to study the orbital stability.\\
\noindent{\it (A)\ For some $p\geq 4$, $k > 0$ and some $z_0$,
Cauchy problem (\ref{eq41}) with initial value $z(x,0) = z_0$ has
a solution
$z(x,t)$ in $[0,T^*)$ with $T^* < \infty$ or $T^* = \infty$.}\\

Denote $\mathcal{S}_\lambda=\{u\in \mathcal{M};\ I(u)=m\}$. Then
with the help of Lemma \ref{le43} and Proposition \ref{pr45}, we
can check in details the proof of Theorem \ref{th33} and get
directly the following orbital stability result.
\begin{Proposition}\label{pr47}
Suppose ${\it (A)}$ holds and $4\leq p < 6$, $\theta>0$. Then for
any $\varepsilon > 0$, there exists $\delta > 0$ such that if
\begin{equation}\label{eq416}
\inf_{u\in\mathcal{S}_\lambda}\|z_0 - u\| < \delta,
\end{equation} then the solution
$z(x,t)$ of (\ref{eq41}) with initial value $z(x,0)=z_0$ satisfies
\begin{equation}\label{eq417}
\inf_{u\in\mathcal{S}_\lambda}\|z(\cdot,t) - u\| < \varepsilon,
\end{equation} for any $t\in[0,T^*)$.
\end{Proposition}
We are now in a position to state the following stronger orbital
stability result.
\begin{Theorem}\label{th48}
Suppose ${\it (A)}$ holds and $4\leq p < 6$, $\theta>0$. Then for
any $\varepsilon > 0$, there exists $\delta > 0$ such that if
\begin{equation}\label{eq418}
\inf_{\eta,\xi\in\R}\|z_0 -  e^{i\eta}u_0(\cdot+\xi)\| < \delta,
\end{equation} then the solution
$z(x,t)$ of (\ref{eq41}) with initial value $z(x,0)=z_0$ satisfies
\begin{equation}\label{eq419}
\inf_{\eta,\xi\in\R}\|z(\cdot,t) - e^{i\eta}u_0(\cdot+\xi)\| <
\varepsilon,
\end{equation} for any $t\in[0,T^*)$.
\end{Theorem}
\pf  It suffices to show that
$\mathcal{S}_\lambda=\{e^{i\eta}u_0(x+\xi)\}$. If
$u\in\mathcal{S}_\lambda$, then
$e^{i\eta}u(x)\in\mathcal{S}_\lambda$ and an argument presented in
the proof of Theorem \ref{th46} implies that $u$ is a solution of
(\ref{eq42}) with some $\mu$. Note that $u(x)$ is a positive
solution of (\ref{eq42}) is equivalent to
$w(x)=\theta^{1\over{p-2}}\mu^{1\over{2-p}}u(\mu^{-{1\over 2}}x)$
is a solution of
\begin{equation}\label{eq420}
-w'' + w - k\theta^{2\over{2-p}}\mu^{1\over{p-2}}(w^2)''w =
w^{p-1},\quad w
> 0,\quad w\in H^1(\R).
\end{equation}
By a result of Ambrosetti-Wang \cite[Page 61]{aw}, we know that
(\ref{eq415}) has a unique positive solution $w_0(x)$ up to a
translation. It follows that
$u_0(x)=\theta^{1\over{p-2}}\mu^{1\over{p-2}}w_0(\mu^{1\over 2}x)$
is the unique positive solution of (\ref{eq415}) up to a
translation. The conclusion follows.$\Box$\\

\end{document}